%% file: FiberTuning.tex
\documentstyle[twoside,12pt,toc_small,bib_toc,psfig]{article}

\input{LaTeX-top}

\newfont{\script}{eusm10 at 12pt}
\def\M{{\mbox{\script M}}}  

\edef\theta{\vartheta}
\long\def\longhide#1{}

\def\newsection#1{
   \section{#1}
}

\newlength\captionwidth
\def\LemImpressFiber{Lemma~2.5}
\def\PropFiberLocConn{Proposition~2.9}
\def\PropLocConnFiber{Proposition~2.10}
\def\PropLocConnJuliaFiber{Proposition~3.6}
\def\CorDisconnectJulia{Corollary~3.7}

\def\ThmPeriodicFiber{Theorem~3.5}
\def\ThmLindeloef{Theorem~A.5}
\def\LemDisconnect{Lemma~A.8}

\def\ThmBranch{Theorem~2.2}
\def\PropFiberLocConn{Proposition~3.5}
\def\ThmHypCompFiberTrivial{Theorem~5.2}
\def\CorMisiuDisconnect{Corollary~6.5}
\def\CorNoNonHypRatRay{Corollary~7.2}
\def\PropPeriodicParaRays{Proposition~7.6}

\title{\vspace{-15mm} On Fibers and Renormalization\\
of Julia Sets and Multibrot Sets}
\author{Dierk Schleicher\\
Technische Universit\"at M\"unchen}
\date{}

\begin{document}

\maketitle
\thispagestyle{empty}
\input{imsmark.tex}
\def\IMSmarkvadjust{-30pt}
\SBIMSMark{1998/13b}{December 1998}{}


\begin{center}
\begin{minipage}{110mm}
\def\toc_vspace{1.5em}
\def\tocname{Contents}
\tableofcontents{0.3em}
\end{minipage}
\end{center}
\vspace{5mm}

\begin{abstract}  
We continue the description of Mandelbrot and Multibrot sets and of
Julia sets in terms of {\em fibers} which was begun in \cite{Fibers}
and \cite{FiberMulti}. The question of local connectivity of these
sets is discussed in terms of {\em fibers} and becomes the question
of triviality of fibers. In this paper, the focus is on the behavior
of fibers under renormalization and other surgery procedures. We show
that triviality of fibers of Mandelbrot and Multibrot sets is
preserved under tuning maps and other (partial) homeomorphisms.
Similarly, we show for unicritical polynomials that triviality of
fibers of Julia sets is preserved under renormalization and other
surgery procedures, such as the Branner-Douady homeomorphisms. We
conclude with various applications about quadratic polynomials and
its parameter space: we identify embedded paths within the Mandelbrot
set, and we show that Petersen's theorem about quadratic Julia sets
with Siegel disks of bounded type generalizes from period one to
arbitrary periods so that they all have trivial fibers and are thus
locally connected.
\end{abstract}

\newpage

\newsection{Introduction}
\label{SecIntro}

This is a continuation of the papers \cite{Fibers} and
\cite{FiberMulti} in which we have introduced the concept fibers for
Multibrot sets, for filled-in Julia sets and, more generally, for
compact connected and full subsets of the complex plane. The idea of
fibers is to use pairs of external rays landing at common points to
cut such a set $K\subset\C$ as finely as possible into subsets (if
the set $K$ has interior points, then extra separations are needed to
cut interior components apart). We call these subsets {\em fibers};
fibers will be defined precisely in Sections~\ref{SecFiberTuning}
respectively \ref{SecJuliaRenorm}. A fiber is called {\em trivial} if
it consists of a single point. In this case, the set $K$ is locally
connected at this point, but triviality of the fiber is a slightly
stronger property. It is equivalent to ``shrinking of puzzle pieces''
(in the sense of Branner, Hubbard and Yoccoz) but independent of any
particular choice of puzzles. It turns out that many proofs of local
connectivity in holomorphic dynamics actually prove triviality of
fibers. For Multibrot sets, there are new direct proofs of triviality
of fibers and thus of local connectivity for all boundary points of
hyperbolic components and for all Misiurewicz points
\cite{FiberMulti}. The investigation of fibers for Julia sets allows,
among other things, to prove that the ``impressions'' of periodic
rays are in many cases only points \cite{Fibers}.  We will build on
these results and discuss how fibers behave under tuning and
renormalization. 

In Section~\ref{SecMultibrot}, we review the definitions of
Mandelbrot and Multibrot sets and discuss renormalization of Julia
sets and tuning maps within Mandelbrot and Multibrot sets. In
Section~\ref{SecFiberTuning}, we show that tuning maps preserve
triviality of fibers: the fiber of any point in a Multibrot set is
trivial if and only if the fiber of the corresponding point in a
``little Multibrot set'' is trivial. Any proof of local connectivity
via shrinking of puzzle pieces for non-renormalizable parameters
in Multibrot sets thus turns automatically into a proof for finitely
renormalizable parameters. A similar result holds for Julia sets,
but more is true: in Section~\ref{SecJuliaRenorm}, we show that
whenever any renormalization of a polynomial with a single critical
point has a locally connected Julia set, then the entire Julia set is
locally connected (and conversely). For example, Petersen \cite{Pt}
has shown that Julia sets of quadratic polynomials with ``bounded
type Siegel disks'' of period one are locally connected; it follows
that this is true for arbitrary periods. This result is discussed in
Section~\ref{SecApplications} together with further results which are
known only for quadratic polynomials; among them is a discussion of
how close the Mandelbrot set comes to being arcwise connected.

The definition of fibers depends of course on the collection of
external rays used to cut the complex plane. For arbitrary compact
connected and full set in $\C$, fibers may behave rather badly and 
there is no universal best choice of external rays. In the beginnings
of Sections~\ref{SecFiberTuning} and \ref{SecJuliaRenorm}, we will
give the exact definitions of fibers for Multibrot sets respectively
for Julia sets together with a review of their most important
properties. These two sections are independent of each other, while
Section~\ref{SecApplications} builds on both of them.

{\sc Acknowledgements}. As its continuation, this paper owes a lot
to the same people which are already mentioned in the
acknowledgements of \cite{Fibers}. I would like to gratefully
emphasize in particular the contribution and continuing support of
Misha Lyubich, John Milnor and the Institute of Mathematical Sciences
at Stony Brook.

\newsection{Tuning and Renormalization}
\label{SecMultibrot}

In this paper, we will be concerned with polynomials having only a
single critical point of possibly higher multiplicity. Milnor has
suggested to call these polynomials {\em unicritical}. They can have
any degree $d\geq 2$, and up to normalization (by affine conjugation)
they can always be written $z\mapsto z^d+c$. We will always suppose
that our Julia sets are normalized in this way. The number $c$ is a
complex parameter which is uniquely determined up to multiplication
by a $d-1$-st root of unity. The filled-in Julia set of such a
polynomial always has exactly $d$-fold rotation symmetry (except if
$c=0$) and is connected if and only if it contains the only critical
point.

The {\em Multibrot set}\/ $\M_d$ of degree $d\geq 2$ is the set of
parameters $c$ for which the Julia set of $z^d+c$ is connected. It
is itself connected and has $d-1$-fold rotation symmetry (see 
\cite{Intelligencer} for pictures of various of these sets). The
case $d=2$ is the familiar Mandelbrot set. Similarly as
\cite{FiberMulti}, the present paper can be read with the case $d=2$
in mind throughout, but the higher degrees do not introduce any new
difficulties. Certain results are known only in the quadratic case;
they are collected in the final section of this paper.

Some of the results in this paper have been floating around the
holomorphic dynamics community in some related form, sometimes
without written proofs. We have included them here to provide
references and proofs, as well as to generalize some of them from the
quadratic case to the case of unicritical polynomials of arbitrary
degrees.

We will make essential use of the theory of {\em polynomial-like
maps} of Douady and Hubbard \cite{Polylike}. A polynomial-like map
is a proper holomorphic map $f\colon U\to V$, where $U$ and $V$ are
open and simply connected subsets of $\C$ such that $\ovl U\subset
V$. Such a map always has a positive degree. Douady and Hubbard
require their maps to have degrees at least two, and this is also the
only case of interest to us. The {\em filled-in Julia set} of a
polynomial-like map is the set of points $z\in U$ which remain in $U$
forever under iteration of $f$; it is always compact, and it is
connected if and only if it contains all the critical points of $f$.
The {\em Julia set} of $f$ is the boundary of the filled-in Julia set.

The Straightening Theorem \cite[Theorem 1] {Polylike} states that
every polynomial-like map $f\colon U\to V$ of degree $d$ is {\em
hybrid equivalent} to a polynomial $p$ of degree $d$ restricted to
an appropriate subset $U'$ of $\C$ such that $V':=p(U')$ contains
$\ovl U'$: that is, there exists a quasiconformal homeomorphism
$\phi \colon V\to V'$ conjugating $f$ to $p$, i.e. $p\circ \phi =
\phi\circ f$ on $U$; this homeomorphism maps the filled-in Julia set
of $f$ onto the filled-in Julia set of $p$, and the complex
dilatation of $\phi$ vanishes on this filled-in Julia set. When the
filled-in Julia set of $f$ is connected, then the polynomial $p$ is
unique up to affine conjugation.

\begin{definition}[Renormalization]
\label{DefRenorm} \lineclear
A unicritical polynomial $p$ of degree $d$ is {\em
$n$-renormalizable} if there exists a neighborhood $U$ of the
critical value such that the holomorphic map $p^{\circ n}\colon U\to
p^{\circ n}(U) $ is polynomial-like of degree $d$ and has connected
Julia set. Any polynomial which is hybrid equivalent to this
polynomial-like map is called a {\em straightened polynomial} of
this renormalization. Such a polynomial is again unicritical of
degree $d$, and in the normalization $z^d+c$ it is unique up to
multiplication of the parameter $c$ by a $d-1$-st root of unity.

The {\em little Julia set} of this renormalization is the filled-in
Julia set of the polynomial-like map. The renormalization is {\em
simple} if the little Julia set does not disconnect any of its images
under $p,p^{\circ 2}, \ldots, p^{\circ(n-1)}$.
\end{definition}
\remark
McMullen, in his general investigation of renormalization of
quadratic polynomials \cite{McBook1}, distinguishes three types of
renormalization: he says that a renormalization of period $n$ with
little Julia set
$K'$ is of
\begin{itemize}
\item{\bf disjoint type}
if $K'$ is disjoint from $p(K'), p^{\circ 2}(K'), \ldots,
p^{\circ(n-1)}(K')$;
\item{\bf $\beta$-type}
if $K'$ intersects any of $p(K'), p^{\circ 2}(K'), \ldots,
p^{\circ(n-1)}(K')$ only at a point which corresponds to the landing
point of a fixed ray of the straightened polynomial;
\item{\bf $\alpha$-type}
if $K'$ intersects any of $p(K'), p^{\circ 2}(K'), \ldots,
p^{\circ(n-1)}(K')$ only at a point which corresponds to the unique
fixed point of the straightened polynomial that is not the landing
point of a fixed ray.
\end{itemize}
McMullen shows that any renormalization of a quadratic polynomial
has one of these three types. A renormalization is simple if and
only if it is of one of the first two types. The third type is also
known as {\em crossed renormalization}, cf.\ \cite{CrossRenorm}.

A quadratic polynomial with connected Julia set has two fixed points:
one is the landing point of the unique fixed ray at angle $0$; this
fixed point is called the $\beta$-point. The other fixed point is
called the $\alpha$-point; thus the names of the renormalization
types. These fixed points coincide exactly when the polynomial is
conjugate to $z^2+1/4$.

Any unicritical polynomial of degree $d\geq 2$ has $d$ finite fixed
points (counting multiplicities). If the Julia set is connected,
then the $d-1$ fixed rays land at $d-1$ distinct fixed points, which
are the analogues of the $\beta$-fixed points: if two rays were to
land together, they would cut $\C$ into two parts, one of which does
not contain the critical point of the polynomial. This part would
then have to map homeomorphically onto itself, which is a
contradiction. The remaining fixed point is called $\alpha$; it may
be attracting, indifferent, or repelling, and it may coalesce with
one of the $d-1$ other fixed points (this happens exactly at the
$d-1$ cusps of the unique hyperbolic component of period $1$). There
is thus at most one fixed point which is not the landing point of a
fixed ray. Statement and proofs of McMullen's classification now
generalize to unicritical polynomials of arbitrary degrees. 

It turns out that the sets of simply $n$-renormalizable unicritical
polynomials are organized in the form of little copies of Multibrot
sets: disjoint types occur at {\em primitive copies} and
$\beta$-types occur at {\em non-primitive} or {\em satellite} copies
of Multibrot sets (see below). In order to describe this, we need to
introduce tuning maps of Multibrot sets.

\begin{definition}[Tuning Map]
\label{DefTuning} \lineclear
A {\em tuning map} of period $n$ of a Multibrot set $\M_d$ is a
homeomorphism $\Psi:\M_d\to \Psi(\M_d)\subset\M_d$ such that, for
every $c\in \M_d$, the polynomial $z\mapsto z^d+\Psi(c)$ is
$n$-renormalizable, and the corresponding polynomial-like map is
hybrid equivalent to the polynomial $z\mapsto z^d+c$. 
\end{definition}

Douady and Hubbard have shown the following theorem in the case of
the Mandelbrot set, but they have not published a complete proof.
It can be found in Section~10 of the recent thesis \cite{Hai} of
Ha\"{\i}ssinsky. Underlying is the theory of ``Mandelbrot-like
families'' \cite{Polylike} which is in general known to work only for
quadratic polynomials. However, the difficulties do not lie in the
degree of the maps but in the number of independent critical points.
Therefore, their theory applies also to families of unicritical
polynomials of arbitrary degrees.

\begin{theorem}[Tuned Copies of Multibrot Sets]
\label{ThmTuningSatellite} \lineclear
For every hyperbolic component of period $n$ of any Multibrot set\/
$\M_d$, there exists a tuning map of period $n$ sending the unique
component of period\/ $1$ to the specified hyperbolic component. This
tuning map is unique up to precomposition with a rotation by
$k/(d-1)$ of a full turn, for some\/ $k\in\{0,1,2,\ldots, d-1\}$. 

If the component is primitive, then the tuning map can be extended
as a homeomorphism to a neighborhood of\/ $\M_d$ onto its image. If
the component is non-primitive, then the tuning map can be extended
as a homeomorphism to a neighborhood of\/ $\M_d$ minus one of the
$d-1$ roots of the period one component. In both cases, the extension
is  no longer unique.
\end{theorem}

\hide{
In order to describe the boundary of little Julia set in a
renormalizable Julia set, we need the following result.
\begin{proposition}[Boundary of Little Julia set]
\label{PropLittleJuliaBdy} \lineclear
Suppose that some unicritical polynomial $p$ with filled-in Julia set
$K$ is $n$-renormalizable, such that there exist open simply
connected domains $U,V\subset\C$ for which the restriction $p^{\circ
n}\colon U\to V$ is polynomial-like of the same degree $d$ with
connected little Julia set. Then there is a finite collection of
periodic dynamic ray pairs landing at points in $U$ and bounding an
open subset $U'\subset U$ with the property that the little Julia set
of the renormalization is exactly the set of points in $K\cap \ovl
U'$ which never leave $\ovl U'$ under iteration of $p^{\circ n}$.
\end{proposition}
\proof
}

\newsection{Fibers of Multibrot Sets and Tuning}
\label{SecFiberTuning}

The main result in this section is that triviality of fibers of 
Multibrot sets is preserved under tuning. Julia sets are discussed in
the next section, and further similar statements for the Mandelbrot
set will be given in Section~\ref{SecApplications}.

We begin by defining fibers of Multibrot sets. What follows is a
brief review of the definition in \cite{Fibers}, simplified using
results from \cite[Section~7]{FiberMulti}. We will use parameter rays
at periodic external angles; they are known to land at parabolic
parameters. Two parameter rays at periodic angles are called a {\em
ray pair} if they land at a common point. We say that two points in
the Multibrot set are {\em separated} by this ray pair if they are
different from the landing point and if they are in different
connected components of $\C$ minus ray pair and landing point. In
order to be able to separate points within hyperbolic components, we
also allow separations by two parameter rays at periodic angles which
land at the boundary of a common hyperbolic component, together with
a simple curve within this hyperbolic component connecting the two
landing points. A {\em fiber} of a Multibrot set is an equivalence
class of points which cannot be separated from each other. A fiber is
called {\em trivial} if it consists of a single point.

This definition might seem special on two accounts: why are only
hyperbolic components allowed for separation lines, excluding
non-hyperbolic components, and why use only parameter rays at
periodic angles, excluding preperiodic angles? In \cite{FiberMulti},
fibers have indeed been defined using parameter rays at periodic and
preperiodic angles, and separation lines through arbitrary interior
components were allowed. It turns out, however, that parameter rays
at rational external angles never land at non-hyperbolic components
by \cite[\CorNoNonHypRatRay]{FiberMulti}, so allowing them would
change nothing except possibly confuse. As to parameter rays at
preperiodic angles, they are simply not necessary: omitting them does
not change fibers at all \cite[\PropPeriodicParaRays]{FiberMulti}.
Since preperiodic parameter rays sometimes need special attention and
separate arguments, it is simply a matter of convenience to exclude
them.

Fibers of Multibrot sets are known to be compact, connected and full.
If any fiber is trivial, then the Multibrot set is locally connected
at this point; this is one of the fundamental properties which make
fibers useful: see \cite[\PropFiberLocConn]{Fibers} or
\cite[\PropFiberLocConn]{FiberMulti}. All points on closures of
hyperbolic components, as well as all Misiurewicz points, have
trivial fibers. This includes all the landing points of parameter
rays at rational angles. (For the Mandelbrot set, many more points
are known to have trivial fibers; see Yoccoz'
Theorem~\ref{ThmYoccozMandelFiber}.)

For the purposes of this paper, there is another possible
simplification: we can ignore separations through hyperbolic
components and just restrict to periodic parameter ray pairs. Of
course, this no longer allows to separate hyperbolic components, but
all fibers on closures of hyperbolic components are trivial and we no
longer have to discuss them. Any parameter which is not on the
closure of a hyperbolic component can be separated from any
hyperbolic component by a periodic ray pair, so its fiber can simply
be constructed using periodic parameter ray pairs only. This will
simplify some discussions below.

It is well known but not so well published that simply renormalizable
parameters of any particular type in the Mandelbrot set are organized
in the form of little Mandelbrot sets (compare Douady and
Hubbard~\cite{Polylike}, Douady~\cite{DoAngles},
McMullen~\cite{McBook1}, Milnor~\cite{MiOrbits}, Riedl~\cite{Riedl}
or Ha\"{\i}ssinsky~\cite{Hai}). Each little Mandelbrot set comes with
a parameter ray pair separating this little Mandelbrot set from the
origin, and the angles of this ray pair are periodic with the same
period $n$ as the renormalization period. In fact, any little
Mandelbrot set can be separated from the rest of $\M_2$ by this
periodic parameter ray pair and a countable collection of parameter
ray pairs. The same statement is true for all the Multibrot sets. 
Little Multibrot sets within $\M_d$ have decorations attached only at
tuned images of Misiurewicz points at external angles $a/d^n$ for
integers $a$ and $n$. Here is a precise statement of the result we
need. We will not include a proof here, although there is no precise
classical reference; compare the references cited above. 

\begin{theorem}[Decorations at Little Multibrot Sets]
\label{ThmLittleMandelBoundary} \lineclear
For any little Multibrot set\/ $\M'$ set within the entire Multibrot
set\/ $\M_d$, there is a countable set $B$ of parameters with the
following property: for any parameter $c\in \M_d-\M'$, there is a
rational parameter ray pair landing at a point $c'\in B$ and
separating $c$ from $\M'-\{c'\}$. This set $B$ consists of a single
parabolic parameter and countably many Misiurewicz points. 

Under the tuning map $\tau:\M_d\to\M'$, the parabolic point in $B$ is
the image of the landing point of the parameter ray at angle
$k/(d-1)$ for some integer $k$, and the Misiurewicz points in $B$ are
exactly the images of the landing points of the parameter rays at
angles $a/d^n$, for any pair of positive integers $a$ and $n$. No
point in $\tau^{-1}(B)$ disconnects $\M_d$.
\qedd
\end{theorem}

\begin{corollary}[Fiber Triviality of $\M_d$ Preserved Under Tuning]
\label{CorFiberTuning} \lineclear
The fiber of any point in the Multibrot set is trivial if and only if
the fiber of the image point under any and all of the tuning maps is
trivial.
\end{corollary}
\proof
With the right background about tuning as given in the previous
theorem, this proof is quite easy. Let $\tau\colon\M_d\to\M'$ be a
tuning homeomorphism and consider a point $c\in\M_d$ together with
$c':=\tau(c)\in\M'$. By \cite[Sections~5 and 6]{FiberMulti}, the
landing points of all the rational parameter rays of $\M_d$ have
trivial fibers. 

First observe that the fiber of any point in $\M'$ must entirely be
contained in $\M'$: by Theorem~\ref{ThmLittleMandelBoundary}, any
point in $\M'$ is either the landing point of a rational parameter
ray and its fiber is trivial, or it can be separated from any other
point in $\M_d-\M'$ by a periodic parameter ray pair. Therefore, we
have to show that $c$ can be separated from any point in $\M_d-\{c\}$
if and only if $c'$ can be separated from any point in $\M'-\{c'\}$. 
By the discussion at the beginning of this section, we only have to
consider separations by periodic parameter ray pairs.

The main idea is to transfer separating ray pairs in $\M_d$ to
separating ray pairs in $\M'$ and back. These ray pairs consist of
parameter rays at periodic angles which land at parabolic parameters,
and all their landing points have trivial fibers. A parabolic
parameter remains parabolic after tuning, and the number of parameter
rays landing there will always be one or two (one at  co-roots, two
at the root of any hyperbolic component of period at least two).
Roots or co-roots map under tuning to roots respectively co-roots. It
follows that periodic ray pairs preserve their separation properties:
when a periodic ray pair separates $c_1$ and $c_2$, then the ray pair
landing at the tuned image of the landing point will separate
$\tau(c_1)$ and $\tau(c_2)$ in $\M'$, and conversely. 
\hide{The same is
true for separation lines through interior components, which
necessarily use parameter rays at periodic angles. 
}

\hide{
For Misiurewicz points, there is a slight complication, as the number
of parameter rays might increase under tuning. Still, the only thing
we need is that whenever a Misiurewicz point separates $\M_d$ into
some number of connected components, then we have just as many
preperiodic parameter rays landing there, separating all these
connected components; this is {\CorMisiuDisconnect} in
\cite{FiberMulti}. After tuning, the image point will separate $\M'$
into just as many connected components. They must all be in different
connected components of $\M_d$ minus the Misiurewicz point because
otherwise the Misiurewicz point would be on the boundary of an
interior component of $\M_d$. This is impossible because boundary
points of hyperbolic components have indifferent orbits, which is not
the case for Misiurewicz points, and triviality of fibers of
Misiurewicz points implies that they cannot be on the boundary of a
non-hyperbolic component. Again by {\CorMisiuDisconnect} in
\cite{FiberMulti}, we will then have enough separating parameter rays
at the tuning image of the Misiurewicz point. (The image of the
Misiurewicz point might separate $\M_d$ into more connected
components, and we have accordingly more parameter rays landing
there.)
}

We conclude the following: whatever separations are possible in
$\M_d$, we will obtain the corresponding separations in $\M'$. The
converse statement works in the same way. This is what we needed to
prove.
\qed

\remark
The nice thing about this proof is that, while proving that the
fibers of arbitrary points in $\M_d$ are trivial, we only have to
worry about periodic parameter rays and their behavior under tuning,
and these do not cause any difficulty. In a sense, we thus do not
prove that new fibers are trivial, but we just show that triviality
is ``preserved'' under certain maps. For example, it is sometimes
possible to conclude certain results for real parameters, using
purely real methods (compare the final section). We can then carry
over these results to many non-real parameters.

\hide{
The original definition of fibers in \cite{FiberMulti} used
separations not only via periodic parameter ray pairs but also via
preperiodic parameter ray pairs and via separation lines through
hyperbolic components. These extra separations do not make fibers
smaller except at hyperbolic components where fibers are known to be
trivial anyway. However, this proof applies even in this extended
definition of separation lines.
}

The very same proof shows that the fiber of any point within some
Julia set is trivial if and only if the fiber of the corresponding
point within a renormalization of the Julia set is trivial. But for
unicritical polynomials, much more is true: if the fibers of all the
points within the ``little Julia set'' are trivial, then all the
fibers within the entire Julia set are trivial. This will be the main
conclusion in the next section.

\newsection{Fibers of Julia Sets and Renormalization}
\label{SecJuliaRenorm}

In this section, we will discuss fibers of connected filled-in Julia
sets of unicritical polynomials, always assuming them to be monic. To
define fibers, it is easiest to first consider only a filled-in Julia
set without interior. Let $Q$ be the set of dynamic rays at rational
angles. All these rays are known to land at repelling periodic
respectively preperiodic points. We will only consider rays which
land together with at least one more rational ray; a pair of rational
rays landing at a common point is again called a {\em ray pair}. Any
such ray pair cuts $\C$ and the filled-in Julia set into two
(relatively) open sets, and we say that points from different parts
are {\em separated} by this ray pair. No ray pair separates its
landing point from anything. The {\em fiber} of any point $z$ is the
collection of all points in the filled-in Julia set which cannot be
separated from $z$. The fiber is {\em trivial} if it consists of $z$
alone.

If the filled-in Julia set has interior, then we have to extend the
definitions. We will ignore the attracting or parabolic cases because
such Julia sets are known to be locally connected anyway. The
remaining case is that of a Siegel disk. Let $Q$ be the set of
dynamic rays at rational angles, plus the grand orbit (backward orbit
of the entire forward orbit) of all the dynamic rays landing at the
critical value, if any; their external angles are necessarily
irrational. A {\em separation line} is either a pair of rays in $Q$
landing at a common point, or a pair of rays in $Q$ landing at the
boundary of the same Fatou component, together with a simple curve
within this Fatou component connecting the landing points of these
two rays. (Even if two dynamic rays land on the boundary of the same
Fatou component, it is not clear that it must be possible to connect
their landing points by a curve within this Fatou component.) We now
proceed as above: every separation line cuts the complex plane and
the filled-in Julia set into two parts and separates the points on
different sides. The fiber of any point $z$ is then the collection of
all points within the filled-in Julia set that cannot be separated
from $z$, and it is trivial if it is the set $\{z\}$.

We have shown in \cite[Section~3]{Fibers} that the fiber of any point
$z$ in a filled-in Julia set is a compact, connected and full subset
of $\C$ and that the relation ``is in the fiber of'' is symmetric.
The dynamics maps any fiber onto a unique image fiber, either
homeomorphically or as a branched cover according to whether or not
there is a critical point in the fiber. If the fiber of some point
$z$ is trivial, then the filled-in Julia set is locally connected at
this point. In particular, triviality of all fibers implies local
connectivity of the entire set. Conversely, if any unicritical
filled-in Julia set is locally connected, then all its fibers as
constructed above are trivial: the necessary separation lines need
dynamic rays at rational angles and rays on the grand orbit of the
unique ray landing at the critical value.

It is not always true that the relation ``is in the fiber of'' is an
equivalence relation. This is true as soon as all the landing points
of the rays used in the construction have trivial fibers, but this
condition is not always satisfied. For example, it is not true when
there is a Cremer point. The fiber of any repelling periodic or
preperiodic point is trivial as soon as its entire forward orbit can
be separated from the critical value and all periodic Siegel disks
(if any) \cite[\ThmPeriodicFiber]{Fibers}. Important examples are
infinitely renormalizable unicritical polynomials of any degree:
fibers are constructed using rational rays, their landing points
always have trivial fibers, and we have an equivalence relation.

The main goal in this section is to prove that triviality of fibers
is preserved under renormalization: all fibers of any filled-in Julia
set are trivial if and only if all the fibers of any of its
renormalized Julia sets are trivial. The transfer of fiber triviality
is done in two steps, corresponding to three sets: a ``big'' Julia
set $K$, a ``little'' Julia set $K'$ which is an invariant subset,
and a ``model'' Julia set $K_1$ which is hybrid equivalent to the
little Julia set; it is the filled-in Julia set of the renormalized
polynomial. The first step shows that the big Julia set has trivial
fibers whenever the little Julia set does, using separation lines
from the big Julia set. The second step then relates triviality of
fibers of the little Julia set to its model. These three sets are
illustrated in Figure~\ref{FigRenorm}.

\begin{proposition}[Trivial Fibers in Dynamic Subsets]
\label{PropFiberSubsetsA} \lineclear
Let\/ $p$ be a unicritical polynomial with connected filled-in Julia
set\/ $K$. Consider a finite set of dynamic ray pairs at rational
angles landing at repelling orbits and let\/ $U$ be one of the
connected components of\/ $\C$ with these ray pairs removed. Fix any
positive integer $n$. This defines a map\/ $\tilde p\colon K\cap\ovl U
\to \C$ via $\tilde p(z)=p^{\circ n}(z)$. Let\/ $K'$ be the compact
subset of points in $K\cap\ovl U$ which never leave\/ $\ovl U$ under
iteration of the map\/ $\tilde p$ and suppose that\/ $K'$ contains
the critical value. Let\/ $Q$ be a countable forward and backward
invariant set of dynamic rays of\/ $p$ which contains the rays
bounding\/ $U$. Then, for this choice of dynamic rays, $K$ has
trivial fibers if and only if $K'$ does.
\end{proposition}

\begin{figure}[htbp]
\centerline{
\psfig{figure=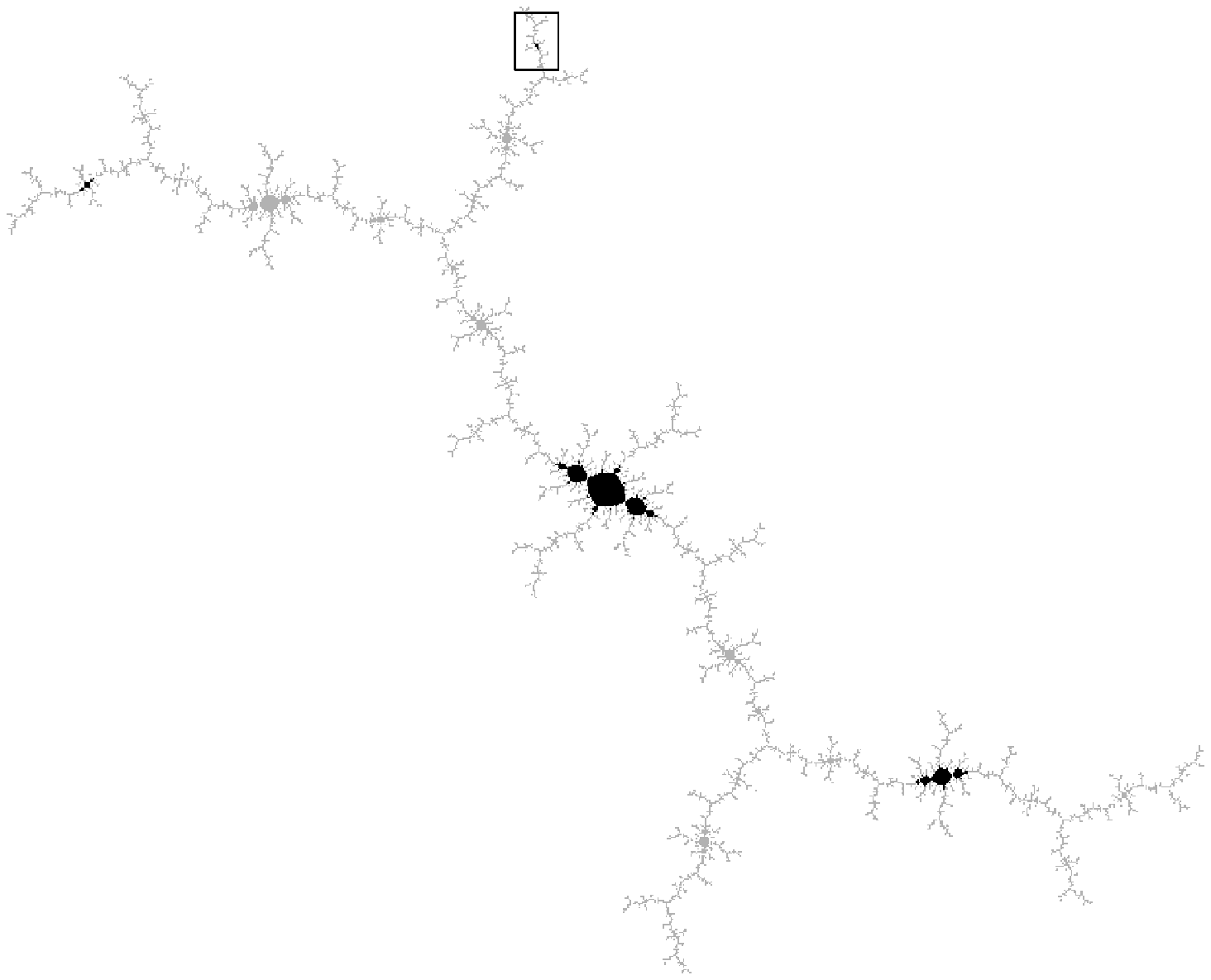,width=110mm}
}
\centerline{
\psfig{figure=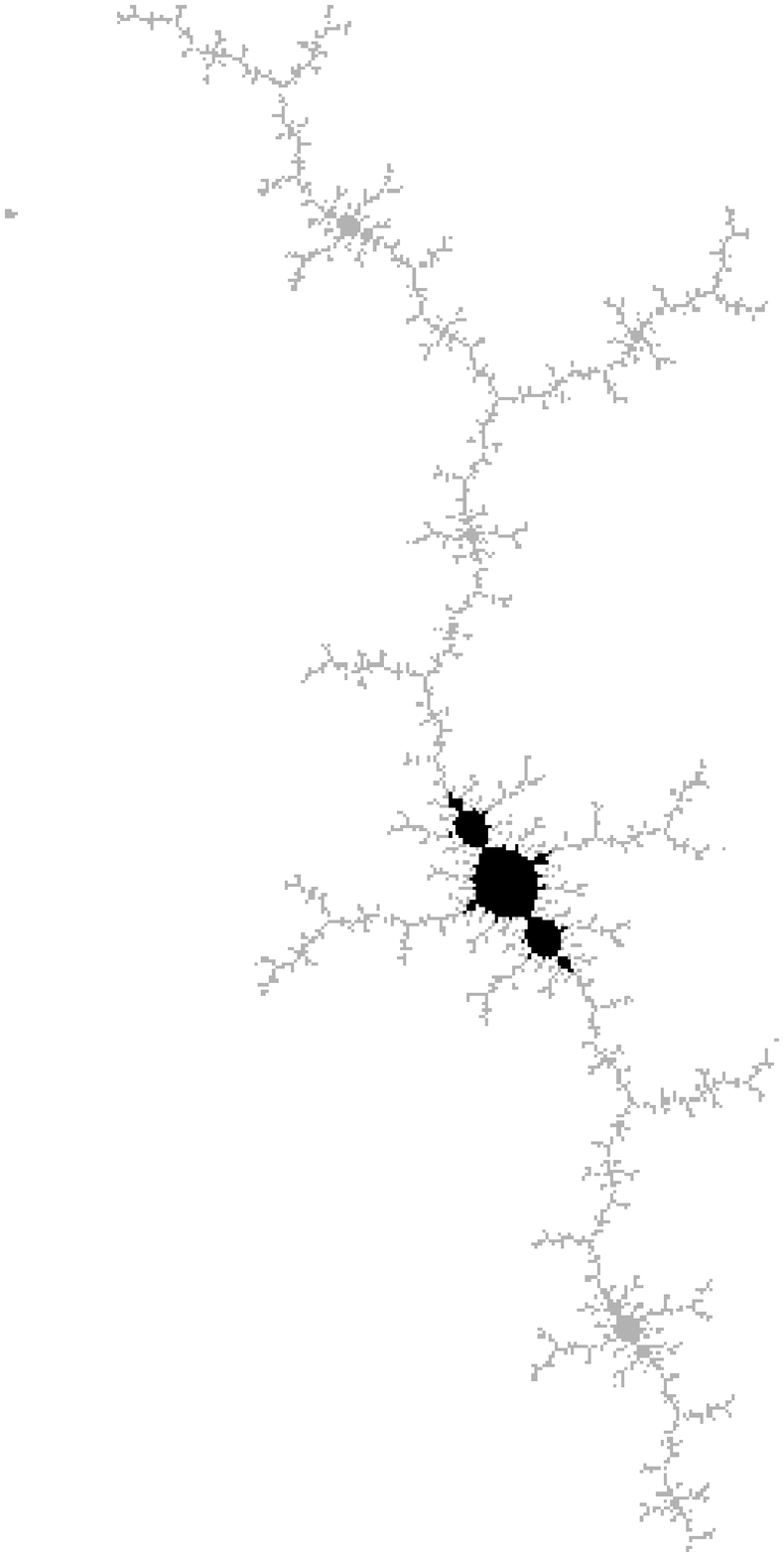,width=40mm}
\psfig{figure=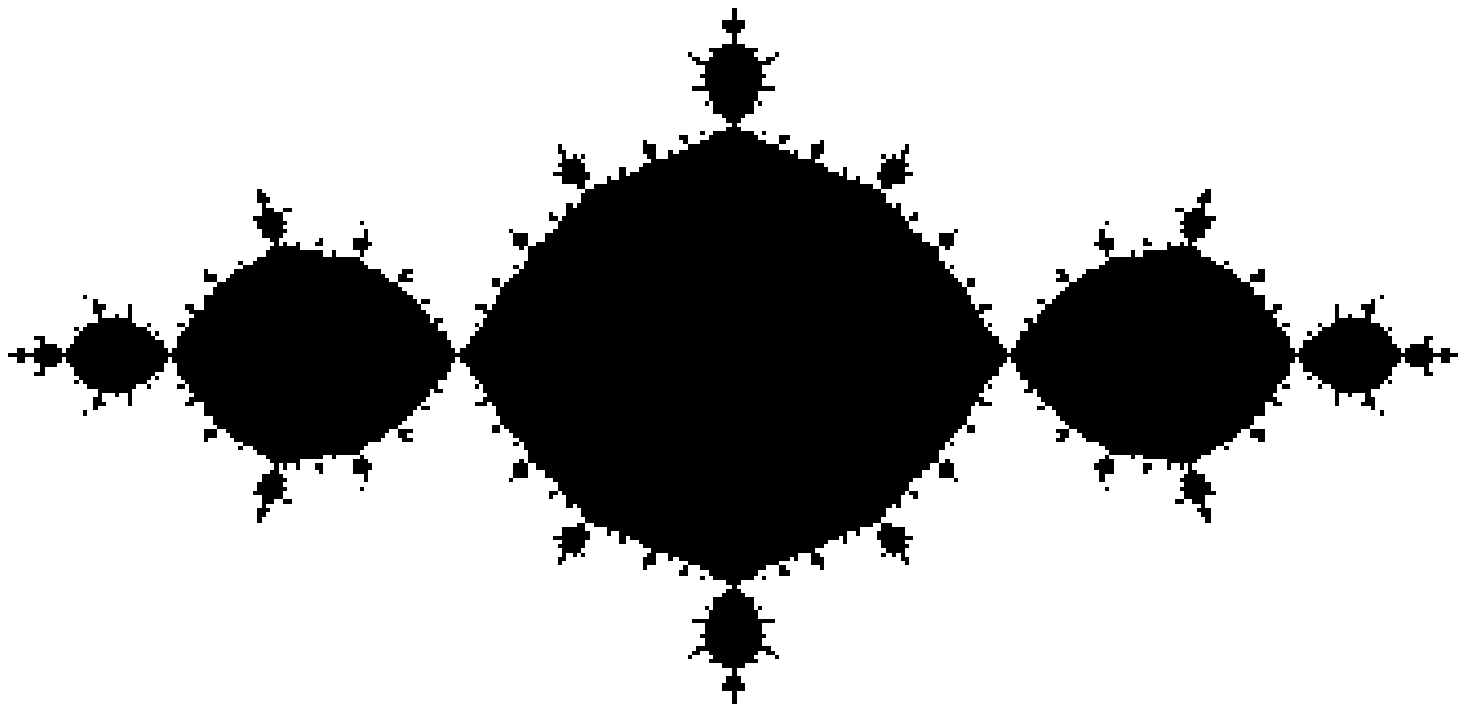,width=70mm}
} 
\LabelCaption{FigRenorm}{Illustration of
Proposition~\protect\ref{PropFiberSubsetsA}. Top: a ``big Julia
set'' $K$ (grey) con\-tai\-ning a ``little Julia set'' $K'$ which
is periodic of period $4$. The little Julia set and its three 
periodic forward images are indicated in black; together, they form
the set $K''$. Bottom left: a blowup near the critical value,
containing the little Julia set $K'$ (black). Lower right:
a model Julia set $K_1$ which is hybrid equivalent to $K'$.}
\end{figure}

\remark
It may be helpful to rephrase the main conclusion of the proposition
as follows: suppose that the dynamic rays (in $Q$) of $K$ supply
enough separation lines so as to disconnect the subset $K'$ in such a
way that all fibers are points. Then the rays in $Q$ disconnect all
of $K$ into fibers consisting of single points. 

The assumption that rays bounding $U$ land on repelling orbits is
largely for convenience. We could as well allow the rays to land on
parabolic orbits. However, since all unicritical polynomials with
parabolic orbits are known to be locally connected and thus to have
trivial fibers anyway, we would not gain much. 
\proof
If $K$ has trivial fibers, then any two of its points can be
separated. In particular, any two points in $K'$ can be separated, so
all the fibers of $K'$ are trivial.

The converse statement requires more work. We assume that $K'$ has
trivial fibers and contains the critical value. Let $S$ be the finite
set of landing points of the rational ray pairs from the statement
and let $S'$ be the union of the points in $S$ with all their forward
orbits, which is still a finite set.

First we assume that the polynomial $p$ has no interior. Then fibers
are constructed using pairs of dynamic rays at rational angles, and
they land at repelling periodic or preperiodic points. Any repelling
periodic point in $K'$ can be separated from the critical value since
the fiber of the critical value is trivial, and any repelling
periodic point outside of $K'$ can be separated from $K'$ and thus
from the critical value as well. Therefore, all the repelling
periodic and preperiodic points have trivial fibers by the result
mentioned above. This applies in particular to all the points in $S$.

Any point $z$ within $K'$ can be separated from any other point
within $K'$ because fibers of $K'$ are trivial. The point $z$ can
also be separated from any point $z'\in K-K'$: if $z\in S$ or if $z$
ever maps into $S$, then its fiber is trivial anyway; otherwise,
there is a finite iterate of $\tilde p$ which sends $z'$ outside of
$\ovl U$, while it obviously leaves $z$ within $K'$, away from the
boundary of $U$. The boundary ray pairs of $U$ thus separate the
corresponding forward images of these two points, so $z$ and $z'$ are
in different fibers. Therefore, all the points in $K'$ have trivial
fibers. Since the dynamics preserves triviality of fibers, all the
points which eventually map into $K'$ have trivial fibers.

Finally, we have to consider the case that $z$ is not in $K'$ and
never maps into $K'$. Let 
\[
K'':=
\hide{
\bigcup_{j=1}^{s}\bigcup_{k=0}^{n_j-1} p^{\circ k}(K'\cap\ovl{U_j})=
}
\bigcup_{k=0}^\infty p^{\circ k}(K')
\,\,.
\]
This set is compact and forward invariant under the dynamics of $p$:
it is the union of all the finitely many forward images of $K'$.
The set $K''$ contains only points which eventually visit $K'$ under
iteration of $p$, so all fibers in $K''$ are trivial. The set
$X:=\C-K''$ is open and carries a unique normalized hyperbolic
metric. Moreover, $X':=p^{-1}(X)$ is a proper subset of $X$, so at
every point in $X'$ the hyperbolic metric with respect to $X'$
strictly exceeds the metric with respect to $X$. The map $p\colon
X'\to X$ is an unbranched covering, hence a local hyperbolic
isometry, so any branch of the inverse is again a local hyperbolic
isometry. With respect to the hyperbolic metric of $X$ in domain and
range, any local branch of $p^{-1}$ contracts the hyperbolic metric
uniformly on compact sets.

Next we claim that the orbit of $z$ has an accumulation point in $X$.
Let 
\[ 
A:= \bigcap_{N>0} \ovl{ \bigcup_{k\geq N} p^{\circ k}(z)}
\]
be the accumulation set of $z$ (i.e. its $\omega$-limit set). Since
$z$ never maps into $K'$ by assumption, the set $A$ is obviously
contained in $\ovl{K-K''}$. If $A$ does not contain a point in $X$,
then the fact $\partial X\cap K\subset \cup_k p^{\circ k}(S)$ implies
that $A$ must consist entirely of points in $S'$, which is a finite
set. But since the orbits of these points are repelling, this is
possible only if the orbit of $z$ falls exactly onto a point in $S'$,
so it lands in $K''$, contrary to our assumption.

Let then $a\in X$ be an accumulation point of the orbit of $z$. Then
for every connected component of $K''$, there is a ray pair in $Q$
separating $a$ from this connected component and avoiding the
finitely many points in $S'$. These finitely many ray pairs will
bound an open neighborhood of $a$, which we will denote by $V$. The
closure of this neighborhood will not meet $K''$, so the hyperbolic
diameter of $V$ in $X$ is finite, hence less than some number
$C<\infty$. For every integer $m>0$, there is an integer $M>m$ such
that the $M$-th iterate of $z$ is in $V$, having mapped at least $m$
times through $V$. Denoting the fiber of $z$ by $Y$, then $f^{\circ
k}(z)\in V$ implies $f^{\circ k}(Y)\subset V$. Therefore, the
hyperbolic diameter of the corresponding forward image of $Y$ is at
most $C$. All the time, the fiber $Y$ is mapped forward
homeomorphically, and the pull-back from forward images of $Y$ back
to $Y$ can only decrease hyperbolic distances. Therefore, the
hyperbolic diameter of $Y$ is at most $C$. Even better: every time
$Y$ maps through $V$, the hyperbolic metrics with respect to $X'$ and
$X$ differ by a definite factor $\alpha<1$, so the hyperbolic
diameter of $Y$ must in fact be less than $\alpha^m C$. Since this is
true for any $m>0$, the fiber $Y$ must be a point, so the fiber of
$z$ is trivial.

It remains to consider the case that the filled-in Julia set has
interior. We only have to consider the case that the polynomial $p$
has a Siegel disk. Since the critical orbit remains in $K''$ and must
accumulate at the boundary of the periodic Siegel disks, all the
periodic Siegel disks are contained in $K''$. The only difference to
the first case lies in the fact that in order to prove that repelling
periodic points have trivial fibers, it must be possible to separate
them from the critical value as well as from all periodic Siegel
disks.

First we show that no boundary of a Siegel disk can contain a
periodic point. Indeed, if there is a repelling point on this
boundary, then it is contained in $K''$ and its fiber is trivial by
assumption. Then  an internal ray of the Siegel disk must land at
this periodic point, but this is impossible: if the ray is periodic,
then the rotation number of the Siegel disk will be rational;
otherwise, more than one internal ray has to land at the same
periodic point of the Siegel disk, and both is a contradiction. 

Since fibers in $K''$ are trivial, every repelling periodic point in
$K''$ can be separated from all the boundaries of the periodic Siegel
disks and from the critical value. This is obviously also true for
repelling periodic points outside of $K''$. Again, all the repelling
periodic and preperiodic points have trivial fibers. The proof now
proceeds exactly as above. 
\qed

The previous proposition was written so that it captures the
construction of renormalization of Julia sets. For other applications
we have in mind, such as the Branner-Douady homeomorphisms, it will
be necessary to phrase the statement somewhat more generally so that
the set $U$ consists of several pieces on which the map $\tilde p$ is
defined as a different iterate of $p$. This requires a compatibility
condition for the dynamics, but the arguments are exactly the same.
In fact, the proof above has been written so that it applies
literally to the following variant of the proposition.

\hide{
\remark
The following proposition is phrased in a more general way than
necessary for tuning because we will also need it for other
applications such as the Branner-Douady homeomorphisms. Statement
and proof are perhaps best read first with the tuning case in mind,
which amounts to $j=1$, $U=U_1$, $\tilde p = p^{\circ n_1}$ (and
without the need for a compatibility condition on the boundary). Then
$K$ is a ``big Julia set'', $K'\subset K$ is a ``little Julia set''
of a renormalization and $K''$ is the entire forward orbit of $K'$.
}

\begin{proposition}[Trivial Fibers in Dynamic Subsets II]
\label{PropFiberSubsets} \lineclear
Let\/ $p$ be a unicritical polynomial with connected filled-in Julia
set\/ $K$. Consider a finite set of dynamic ray pairs at rational
angles landing at repelling orbits and let\/ $U_1,\ldots, U_s$ be
some of the connected components of\/ $\C$ with these ray pairs
removed. Let\/ $U$ be the interior of\/ $\cup \ovl U_j$; it need not
be connected. Further, let\/ $n_1, \ldots, n_s$ be positive integers
such that, whenever\/ $z\in K\cap\partial U_j \cap \partial U_{j'}$,
then\/ $p^{\circ  n_j}(z) = p^{\circ n_{j'}}(z)$. This defines a new
map\/ $\tilde p\colon K\cap\ovl U \to \C$ by setting\/
$\tilde p(z) := p^{\circ n_j}(z)$ whenever\/ $z\in K\cap\ovl U_j$.
Let\/ $K'$ be the compact subset of points in\/ $K\cap\ovl U$ which
never leave\/ $\ovl U$ under iteration of the map\/ $\tilde p$.
Suppose that\/ $K'$ contains the critical value. Let\/ $Q$ be a
countable forward and backward invariant set of dynamic rays of\/ $p$
which  contains the rays bounding the\/ $U_j$. Then, for this choice
of dynamic rays, $K$ has trivial fibers if and only if\/ $K'$ does.
\qedd
\end{proposition}

In both variants of the previous proposition, we have compared
triviality of fibers of Julia sets to triviality of fibers of
``little Julia sets'' within these Julia sets. We will now relate
this to triviality of fibers of renormalized Julia sets.

\begin{theorem}[Trivial Fibers For Little Julia Sets]
\label{ThmFiberLittleJulia} \lineclear
Suppose, under the assumptions of
Propositions~\ref{PropFiberSubsetsA} or \ref{PropFiberSubsets}, that
$\tilde p$ on\/ $K'$ is topologically conjugate to a unisingular
polynomial\/ $p_1$ on its filled-in Julia set $K_1$. Then\/ $K$ has
all its fibers trivial if and only if\/ $K_1$ has the same property,
for appropriate choices of the sets of dynamic rays of the
polynomials $p$ and $p_1$ used to define separation lines for\/ $K$
and\/ $K_1$.
\end{theorem}
\remark
The conjugation assures that the unique critical orbit of $p_1$
(other than $\{\infty\}$) has bounded orbit, so $K_1$ and thus $K'$
are connected. We do not require the conjugation to exist in
neighborhoods of $K'$ and $K_1$, and we do not even require it to
preserve the cyclic order of branch points.

Of course, when $K'$ and $K_1$ are homeomorphic, then they are
simultaneously locally connected. Triviality of fibers is somewhat
stronger; the transfer from $K$ or $K'$ to $K_1$ is still easy. The
reverse direction, showing that trivial fibers of $K_1$ imply trivial
fibers of $K'$, is not immediate: 
\hide{
If all the fibers of $K$ are trivial, then in particular all the
fibers in $K'$ are trivial and $K'$ is locally connected by
\cite[\PropFiberLocConn]{Fibers}. Local connectivity is
preserved under homeomorphisms, so $K_1$ is then also locally
connected. By \cite[\PropLocConnFiber]{Fibers}, $K_1$ has all its
fibers trivial for an appropriate choice of $Q_1$; this choice is
specified in \cite[\PropLocConnJuliaFiber]{Fibers}.
For the reverse transfer, the only obvious statement is that local
connectivity of $K_1$ implies local connectivity of $K'$. The
transfer of triviality of fibers from $K_1$ to $K'$ is not immediate:
}
the additional decorations of $K'$ within $K$ might make it difficult
for dynamic rays to land in such a way as to obtain enough separation
lines. Even if we knew that the topological conjugacy and thus the
homeomorphism extended to a neighborhood, the transfer would not be
obvious: dynamic rays of $K_1$ would then turn into invariant curves
of $K$ avoiding $K'$, but they might run through $K-K'$. (Even hybrid
equivalences between polynomial-like maps can have this problem. An
insignificant problem, easily overcome by Lindel\"of's Theorem
(compare \cite[Theorem~3.5]{Ahl2} or \cite[\ThmLindeloef]{Fibers}),
is that hybrid equivalences do not in general map dynamic rays to
dynamic rays.) What we need is, for every dynamic ray of $K_1$ used
in a separation line, a dynamic ray of $K$ which lands at the correct
point; if there are several rays land the same point in $K_1$, then
we need equally many rays of $K$ (recall that we did not require our
homeomorphism to preserve the order of branches).

As an example, consider the real quadratic polynomial with a
superattracting cycle of period $3$: it is $p(z)=z^2+1.75488\ldots$.
The third iterate is of course hybrid equivalent to the map $z\mapsto
z^2$; the fixed ray at angle $0$ of the latter map turns into a
subset of the real line for $p$ (under a symmetric hybrid
equivalence), so part of the image of the ray is in the filled-in
Julia set of $p$. There are two complex conjugate dynamic rays of $p$
landing at the same point, which are both homotopic to the image of
the $0$-ray with respect to the ``little Julia set''. 

\proof
Suppose that all the fibers of $K$ are trivial for some choice of
dynamic rays. Then in particular all fibers in $K'$ are trivial, and
we want to show that all fibers in $K_1$ are trivial. If the
conjugation between $K'$ and $K_1$ extends to neighborhoods of these
sets, then all the dynamic rays used in separation lines in $K'$
transfer to invariant curves outside of $K_1$ landing at boundary
points of $K_1$. By Lindel\"of's Theorem, there are dynamic rays
landing at the same points through the same accesses. Curves within
bounded Fatou components are of course preserved by the homeomorphism
between $K'$ and $K_1$, so every separation line in $K'$ has a
counterpart in $K_1$. If all fibers of $K'$ are trivial, then all
fibers of $K_1$ are trivial as well.

However, if the conjugation does not extend to neighborhoods of the
filled-in Julia sets, for example if it does not respect the cyclic
order of branch points, then we can still argue as follows:
triviality of all fibers in $K$ triviality of all fibers in $K'$ and
thus local connectivity of $K'$; this turns into local connectivity of
$K_1$. But any locally connected unicritical filled-in Julia set has
only trivial fibers, provided fibers are constructed using all
dynamic rays at rational angles plus the grand orbit of the unique
ray landing at the critical value in case there is a Siegel disk;
compare \cite[\PropLocConnJuliaFiber]{Fibers}. 

The other direction of the theorem is less obvious: we assume that
the fibers of $K_1$ are trivial for some collection $Q_1$ of rays
and show that then an appropriate collection of dynamic rays of $p$
makes all the fibers of $K'$ trivial. We need to show that, for
every dynamic ray in $Q_1$ landing at a point $z_1\in K_1$, there
are dynamic rays of $p$ landing at the corresponding point $z$ in
$K'$ and separating all the connected components of $K'-\{z\}$,
which are equally many in number as the connected components of
$K_1-\{z_1\}$ (these numbers are known to be finite, and all the
rays landing at a common point have the same period). Recall that if
there is no Siegel disk, then we construct fibers using all dynamic
rays at rational external angles; if there is a Siegel disk, then we
have to add the grand orbit of the unique ray landing at the critical
value.

As in the proof of Proposition~\ref{PropFiberSubsetsA}, let $K''$ be
the union of the forward images of $K'$. This is a compact and
forward invariant subset of $K$, and any point in $K''$ will visit
$K'$ after finitely many iterations. 
 
We claim that, for an arbitrary point $z\in K'$, any two connected
components $K_1,K_2$ of $K'-\{z\}$ must be in different connected
components of $K-\{z\}$. Suppose the opposite were true. Since the
filled-in Julia set of $p$ is full, the connected components $K_1$
and $K_2$ must then be connected by a bounded Fatou component of $p$
with $z$ on its boundary. But such a Fatou component must be
contained in $K''$ and even in $K'$, which is a contradiction.

The landing point of every periodic or preperiodic ray is on a
repelling or parabolic orbit. Every periodic Fatou component is
attracting, parabolic or a periodic Siegel disk. In the attracting or
parabolic cases, one of the periodic components contains the critical
value, so that all the periodic Fatou components are contained in
$K''$. Periodic Siegel disks are also contained in $K''$: since $K''$
is forward invariant, it contains either the entire cycle of periodic
Siegel disks or none at all, and in the latter case the boundary of
$U$ would separate the orbit of Siegel disks from the critical orbit,
so the critical orbit could not accumulate on the boundary of the
Siegel disks.

\hide{
This Fatou component cannot correspond to an attracting or parabolic
cycle because it would have to absorb a critical orbit, while the
only critical orbit must be within $K'$. The Fatou component cannot
be a Siegel disk, either: the closure of the critical orbit must
contain the boundary of the Siegel disk; since the critical orbit is
in the closed set $K'$, the Siegel disk must be in $K'$ as well and
$K_1$ and $K_2$ are not in different connected components of
$K'-\{z\}$.
}
 
First consider a periodic angle $\theta_1\in Q_1$. Denote its period
by $m$ and its landing point by $z_1$. Then there is a corresponding
periodic point $z$ in $K'$ which is repelling or parabolic. Any
connected component of $K_1-\{z_1\}$, of which there are finitely
many, corresponds to a unique connected component of $K'-\{z\}$.
All connected components of $K'-\{z\}$ are contained in different
connected components of $K-\{z\}$. The finitely many dynamic rays of
$p$ landing at $z$ separate any two of these connected components
\cite[\LemDisconnect] {Fibers}. Therefore, if any two points in $K_1$
can be separated by a ray pair landing at $z_1$, then the
corresponding points in $K'$ can be separated by a ray pair landing
at $z$. When there are separation lines for $K_1$ consisting of
(pre-)periodic dynamic rays and running through bounded Fatou
components, we get corresponding lines for $K'$.

If there is no Siegel disk, then triviality of fibers of $K_1$
implies triviality of fibers of $K'$ and thus, by
Proposition~\ref{PropFiberSubsetsA}, triviality of all fibers of $K$.
This finishes the proof of the theorem, except if there is a Siegel
disk. In that case, we also have to consider the dynamic ray landing
at the critical value of $p_1$. Let $\theta_1$ be its external angle.
The goal is the same as above: we want to construct a dynamic ray of
$p$ which lands at the critical value of $p$ and which has the same
separation properties for $K'$ as it does for $K_1$. The hard part
this time is to show the existence of an access outside $K$ to the
critical value. Lindel\"of's theorem then supplies a ray landing at
$v$, and it will automatically have the right separation properties
because only one ray can land at the critical value.

Denote the critical values of $p$ and $p_1$ by $v$ and $v_1$,
respectively. We first show that the fiber of $v$ in $K$ is trivial,
albeit using temporarily a different collection of separation lines.
Since all the fibers of $K_1$ are trivial, for every point $z_1\in
K_1$ different from $v_1$ there is a curve in $K_1$ connecting two
repelling periodic points such that this curve, together with two
dynamic rays landing at these points, separates $v_1$ from $z_1$.
This separation line differs from usual separation lines in that its
dynamic rays must land at repelling periodic points, but it is not
required to traverse only a single Fatou component. For every point
$z\in K'$ different from $v$, we now obtain a similar separation line
separating $z$ from $v$. Any point in $K-K'$ can easily be separated
from $v$. Any point in $K-\{v\}$ can thus be separated from $v$ by
such a modified separation line. Since $v$ is on the boundary of $K$,
it is in the impression of the dynamic ray at some angle $\theta$.
This impression cannot cross the modified separation lines just
defined (similarly as in the proof of
\cite[\LemImpressFiber]{Fibers}), so the impression of the dynamic
$\theta$-ray is $\{v\}$ and the ray lands at $v$. At this point,
there might be various rays landing at $v$, which will then be
homotopic in the complement of $K'$. We can transport all these rays
forward and backwards, so that in particular we obtain $d$ dynamic
rays landing at the critical point in a symmetric way (or even a
multiple of $d$ rays). 

The set $K_1$ has trivial fibers for $Q_1=\Q/\Z$, extended by the
grand orbit of the ray landing at the critical value in case there is
a Siegel disk (using separation lines in the usual sense). All these
separation lines of $K_1$ have counterparts in $K$ with the
corresponding separation properties for $K'$. Therefore, all the
separation lines made with dynamic rays landing at repelling and
parabolic (pre-)periodic points of $K'$ and at the grand orbit of the
critical value (in case of a Siegel disk) suffice to make all the
fibers of $K'$ trivial. By the previous proposition, all the fibers
of $K$ are trivial. This finishes the proof of the theorem. (If there
is a Siegel disk, we can now conclude that the critical value is the
landing point of exactly one dynamic ray; see
\cite[\CorDisconnectJulia]{Fibers}.)
\qed

We can now conclude that triviality of fibers of Julia sets is
preserved under renormalization.

\begin{corollary}[Triviality of Fibers Preserved Under
Renormalization]
\label{CorFiberRenormJulia} \lineclear
The Julia set for any polynomial in $\M_d$ has trivial fibers if and
only if any or all of its (simple) renormalizations have the same
property.
\end{corollary}
\proof
This conclusion is immediate once the description of simple
renormalization is right. The usual definition for a polynomial $p\in
\M_d$ to be $n$-renormalizable is that there exist two open subsets
$U\subset V\subset \C$ such that $p^{\circ n}\colon U\to V$ is
polynomial-like of degree $d\geq 2$ with connected Julia set
(Definition~\ref{DefRenorm}).

For quadratic polynomials, it is well known (but admittedly I know of
no precise reference; compare Milnor~\cite{MiOrbits}, 
Lyubich~\cite[Section~2]{LyHypDense} or
Ha\"{\i}ssinsky~\cite[Section~10]{Hai}) that then there are two
rational ray pairs, one periodic and one preperiodic from the same
orbit, which bound a subset $U'\subset U$ with the property that the
little Julia set of the renormalization is exactly the set of points
in $\C$ with bounded orbits which never escape from
$\ovl U'$ under iteration of $p^{\circ n}$. (The set $V$ for the
polynomial-like map in the definition of renormalization above can be
obtained from $U'$ by restricting it to any positive equipotential
and ``thickening'' it slightly near the parts of its boundary which
are formed by dynamic rays and their landing points; the set $U$ is
then obtained from $V$ by pulling it back $n$ times under the
dynamics of $p$.)

With this description of renormalization,
Theorem~\ref{ThmFiberLittleJulia} applies and shows that, whenever
the filled-in Julia set of $p$ has trivial fibers, then the fibers of
any of its renormalizations are also trivial, and conversely.

For unicritical polynomials of higher degrees, the same reasoning
works, but one needs $d-1$ preperiodic ray pairs. 
\qed

\remark
The same statement holds also for crossed renormalization, which has
so far been described only in the quadratic case; see below. Quite
generally, the same ideas might give a statement like this: if an
arbitrary polynomial is renormalizable (possibly using several
regions $U$ and $V$ corresponding to different renormalization
periods) so that these renormalizations absorb all the critical
orbits, then all the fibers of the entire Julia set are trivial if
and only if all the fibers of all the renormalized Julia sets are
trivial.

\newsection{Quadratic Polynomials}
\label{SecApplications}

This section will contain further applications of the concepts
introduced above which are in some sense specific to the case of
quadratic polynomials and to the Mandelbrot set $\M=\M_2$. It seems
quite possible that these statements have analogues for Multibrot
sets of higher degrees, but the underlying theorems are so far known
only for degree two. We will prove that triviality of fibers of both
Julia sets and the Mandelbrot set are preserved under homeomorphisms
like those arising in the context of crossed renormalizations or
those introduced by Branner and Douady, and we will also conclude
that quadratic Julia sets with Siegel disks of bounded type have
trivial fibers. 

\begin{corollary}[Trivial Fibers Preserved under Crossed
Renormalization]
\label{CorFibersCrossedRed} \lineclear
The Julia set of a crossed renormalizable polynomial in $\M$ has
trivial fibers if and only if any of its crossed renormalizations has
the same property.
\end{corollary}
\proof
By \cite{CrossRenorm}, any crossed renormalization which is not of
immediate type (so that the little Julia sets intersect each other at
a fixed point) is itself simply renormalizable in such a way that the
renormalized Julia set is crossed renormalizable of immediate type.
We can therefore restrict attention to the immediate case. For this,
the argument is literally the same as for
Corollary~\ref{CorFiberRenormJulia}: the ray pairs bounding the set
$U'$ are constructed explicitly in \cite[Section~3.1]{CrossRenorm}.
\qed

We will now discuss a homeomorphism which has been introduced by
Branner and Douady \cite{BD}: it maps the $1/2$-limb of the
Mandelbrot set homeomorphically onto a subset of the $1/3$-limb. It
was with this homeomorphism in mind that we added the extra notations
in Proposition~\ref{PropFiberSubsets} and 
Theorem~\ref{ThmFiberLittleJulia}.

\begin{corollary}[Trivial Fibers Preserved Under Branner-Douady Maps]
\label{CorFibersBD} \lineclear
The Branner-Douady homeomorphism from the $1/2$-limb into the
$1/3$-limb of the Mandelbrot set preserves the property of Julia
sets that all fibers are trivial.
\end{corollary}
\proof
As for tuning, this is an immediate conclusion with the right setup.
It is more convenient to look at the inverse map from a subset of the
$1/3$-limb onto the $1/2$-limb. 

For any Julia set from the $1/3$-limb, it is well known that the
dynamic rays at angles $1/7$, $2/7$ and $4/7$ land at a single fixed
point, which is called the $\alpha$-fixed point. Similarly, the
dynamic rays at angles $9/14$, $11/14$ and $1/14$ land at $-\alpha$.
Denote the closures of the regions between the $1/7$-ray and the
$2/7$-ray by $Y_1$, between the $2/7$-ray and the $4/7$-ray by $Y_2$
and let $Z_1:=-Y_1$, $Z_2:=-Y_2$. Finally, let $Y_0$ be the ``central
puzzle piece'' in the complement of all these rays: the region
containing the critical point. Then the critical value is always in
$Y_1$. 

Now the set $U$ consists of the interior of the pieces $Y_0\cup Y_1
\cup Z_1$. More precisely, using the notation from
Theorem~\ref{ThmFiberLittleJulia}, let $U_1:=Y_1$ with $n_1=2$, let
$U_2:=Y_0$ with $n_2=1$ and $U_3:=Z_2$ with $n_3=1$. The dynamics of
$\tilde p$, restricted to the filled-in Julia set, is then as
follows: it sends $U_1\cap K=Y_1\cap K$ homeomorphically onto
$(Y_0\cup Z_1\cup Z_2)\cap K$; $U_3\cap K=Z_2\cap K$ lands
homeomorphically on the same image; and $U_2\cap K=Y_0\cap K$ is
mapped in a two-to-one fashion onto $Y_1\cap K$. Therefore, any point
within $U\cap K$ which eventually escapes from $U$ does so through
$Z_1$. It follows that a point in the $1/3$-limb of the Mandelbrot
set is in the range of the homeomorphism from the $1/2$-limb iff the
critical orbit avoids the piece $Z_1$. 

Let $c_1$ be any parameter in the $1/2$-limb and let $c$ be its image
in the $1/3$-limb under the Branner-Douady homeomorphism. Branner
and Douady then show that the filled-in Julia set of $c_1$ is
topologically conjugate to the subset $K'$ of the filled-in Julia
set of $K$ consisting of points which stay in $\ovl U$ forever.
Therefore, Theorem~\ref{ThmFiberLittleJulia} applies and shows that
the filled-in Julia sets of $c$ and $c_1$ both have trivial fibers
whenever one of them does.
\qed
\remark
By work of Levin and van Strien \cite{LvS}, and of Lyubich and
Yampolsky \cite{LY}, it is known that all the Julia sets on the real
axis of the Mandelbrot set are locally connected and thus have
trivial fibers. 

The construction of Branner and Douady gives, with only notational
changes, homeomorphisms from the $1/2$-limb into any $1/q$-limb of
the main cardioid of the Mandelbrot set. The theorem applies to all
of them. We obtain the following corollary:

\begin{corollary}[Julia Sets in Some Spines of the Mandelbrot Set]
\label{CorJuliaSpineBD} \lineclear
There are topological arcs in the Mandelbrot set connecting the
origin to the landing point of any parameter ray at angle $1/2^n$
and running only along parameters for which their Julia sets have
trivial fibers.
\qedd
\end{corollary}

Using further homeomorphisms, one can extend this result to many
more topological arcs. Such homeomorphisms are still covered by our
arguments. Using results of Yoccoz, as well as the fact that
triviality of fibers is preserved under tuning, one can arrive at
this result in a different way. First we recall (and restate) a
result of Yoccoz.

\hide{
By the same token, triviality of fibers is preserved under other
maps in parameter space and other surgery procedures, for example
those arising in the context of crossed renormalization
\cite{CrossRenorm}: the locus of crossed $n$-renormalizable
parameters within any $p/q$-limb of the Mandelbrot set (with $q$
divisible by $n$) is homeomorphic to the $p/q'$-limb for $q'=q/n$,
and this homeomorphism preserves triviality of fibers. Similarly,
triviality of all fibers of a Julia set is also preserved under
crossed renormalization, similarly as in the case of simple
renormalization above. 
}

\begin{theorem}[Trivial Fibers of the Mandelbrot Set at Yoccoz Points]
\label{ThmYoccozMandelFiber} \lineclear
If a parameter in the Mandelbrot set is is not infinitely
renormalizable, then its fiber is trivial.
\qedd
\end{theorem}
\remark
This theorem is a union of two results by Yoccoz, plus the triviality
that fibers within hyperbolic components are trivial. For the case
that all periodic orbits are repelling, Yoccoz' result is
\cite[Theorem~III]{HY} (only the non-renormalizable case is treated
in detail there; the transfer to the finitely renormalizable case
also follows from our Corollary~\ref{CorFiberTuning}). For the case
of indifferent orbits, Yoccoz' theorem is \cite[Theorem~I.B]{HY}; for
a different proof, see \cite[{\ThmHypCompFiberTrivial}]{FiberMulti}.
The usual statement of Yoccoz' theorems is that the Mandelbrot set is
locally connected at these points. However, he proves the stronger
result that the fibers are trivial by proving shrinking of puzzle
pieces, and we will need triviality of fibers. 

Here is another lemma about the Mandelbrot set which will be needed.
\begin{lemma}[Branch Points on the Real Axis]
\label{LemBranchReal} \lineclear
If a Misiurewicz point in $\M$ is on the real axis, or if it is a
tuned image of a point $c\neq -2$ on the real axis, then it is not a
branch point.
\end{lemma}
\proof
Suppose that there is a Misiurewicz point on the real axis of $\M$
which is a branch point. It is then the landing point of as many
preperiodic parameter rays as there are branches at this point, at
least three (in fact, because of symmetry, this number must be even).
In the dynamic plane of this Misiurewicz point, there must then be at
the same number of preperiodic dynamic rays landing at the critical
value. The critical value is real, and by invariance of the real
axis, there is then a repelling periodic real orbit such that at
least four dynamic rays land at each of its points. This is
impossible for various reasons: the combinatorial rotation number of
this orbit must, by symmetry, be equal to its inverse, so it must be
equal to $1/2$ and the various branches could not be permuted
transitively. But this is always true if there are at least three
rays (compare e.g.\ \cite[Lemma~2.4]{MiOrbits} or
\cite[Lemma~2.4]{ExtRays}). Another reason is that the wake of this
orbit must be bounded by two parameter rays landing at a parabolic
parameter with the same combinatorial rotation number. Again,
symmetry of the real line requires this combinatorial rotation number
to be $1/2$, so the only critical orbit could not visit all the
parabolic basins, which is a contradiction.

If a Misiurewicz point on a tuned image of the real axis is a branch
point, then it must acquire more branches in the process of tuning.
If $c$ is the un-tuned Misiurewicz point, then it cannot disconnect
$\M$ by Theorem~\ref{ThmLittleMandelBoundary}. Therefore, if it is on
the real line, then it is the main antenna tip $c=-2$.
\qed

\begin{theorem}[The Mandelbrot Set is Almost Path Connected]
\label{ThmMandelPath} \lineclear
For every parameter in the Mandelbrot set which has a trivial fiber,
or which is on a tuned image of the real line, there is an arc
within the Mandelbrot set connecting this parameter to the origin.
\end{theorem}
\proof
If the Mandelbrot set was locally connected, the claim would follow
simply by the general fact that compact connected locally connected
metric spaces are pathwise connected \cite[Section~16]{MiIntro}.
There are models for the Mandelbrot set which are locally connected:
one of them is Douady's ``pinched disk model'' (compare
\cite{DoCompacts}). Another one is the quotient of the Mandelbrot set
in which all fibers are collapsed to points (compare
\cite[Section~2]{Fibers} or \cite[Section~7]{FiberMulti}). All
locally connected models of the Mandelbrot set are of course
homeomorphic. Another related locally connected model space is
Penrose's ``abstract abstract Mandelbrot set'' \cite{Pe1},\cite{Pe2},
which is a parameter space of kneading sequences (however, it is not
homeomorphic to the previous spaces because only a subset corresponds
to realized kneading sequences, and many of its points correspond to
more than one point in models of the Mandelbrot set. The necessary
modifications of the proof below are only minor).

All these locally connected model spaces come with natural continuous
projections from the actual Mandelbrot set to these spaces. Let
$\M^*$ be a locally connected model of $\M$ and let
$\pi\colon\M\to\M^*$ be such a projection. The inverse image
$\pi^{-1}(c^*)$ over any point $c^*\in\M^*$ corresponds exactly to a
fiber of the Mandelbrot set. Injectivity at any given point is then
equivalent to this fiber being trivial. 

Let $c_0$ be any point in $\M$ with trivial fiber, or on a tuned
image of the real axis. If $c_0$ is on a tuned image of the real
line, we may replace it by a hyperbolic parameter on the same tuned
image of the real line, so that the fiber of the new $c_0$ is
trivial. It is easily possible to connect the old and new $c_0$ along
the tuned image of the real axis, so nothing is lost by the
assumption that the fiber of $c_0$ is trivial. 

Let $\phi\colon [0,1]\to\M^*$ be an injective continuous map
connecting the origin to $c_0$ within the locally connected model
space $\M^*$. We may suppose that $\phi$ has the following property:
for any little model Mandelbrot set $\M'\subset\M^*$ and two points
$c,c'$ on the tuned image of the real axis of $\M'$ which are on the
image of $\phi([0,1])$, the path connects these two points entirely
along the tuned image of the real axis. (For example, all regular
arcs have this property; see below.) We will show that this map lifts
to a continuous map $\psi\colon[0,1]\to\M$ such that
$\phi=\pi\circ\psi$. 

For any $t\in[0,1]$ such that the fiber of $\phi(t)$ is trivial, the
definition of $\psi(t)=\pi^{-1}(\phi(t))$ is clear. In view of
Yoccoz' theorem above, we only have to consider the case that
$c^*=\phi(t^*)$ is infinitely renormalizable. In fact, it is then
infinitely simply renormalizable \cite{McBook1} so that there is an
infinite sequence of positive integers $n_1<n_2<\ldots$ so that $c^*$
is simply $n_k$-renormalizable and every $n_k$ strictly divides
$n_{k+1}$. Let $\M'_k$ be the nested collection of corresponding
embedded Mandelbrot sets. All of these sets contain $c_0$ in their
wakes, and there are three possibilities for each $\M'_k$:
\begin{enumerate}
\item
the little Mandelbrot set $\M'_k$ may contain $c_0$;
\item
a little Mandelbrot set $\M'_k$ may contain $c_0$ in its wake, so
that the ``main antenna tip'' of $\M'_k$ (the tuned image of the
point $-2$) separates $c_0$ from the origin (i.e. two parameter rays
landing at the antenna tip do the separation);
\item
the little Mandelbrot set $\M'_k$ may contain $c_0$ in its wake, but
not within itself and not so that its main antenna tip separates
$c_0$ from the origin.
\end{enumerate}
Suppose that there is some index $k_0$ so that $\M'_{k_0}$ does not
contain $c_0$.  Then there is a Misiurewicz point $B\in\M'_{k_0}$
which is the landing point of at least two parameter rays at
preperiodic angles which separate $c_0$ from the rest of $\M'_{k_0}$
(Theorem~\ref{ThmLittleMandelBoundary}).  Let $c'_0$ be the center of
any hyperbolic component which is not separated from $c_0$ by the
parameter rays landing at $B$; such a component exists because only
finitely many rays land at $B$ and any periodic parameter ray will
find such a component. As far as the desired path within $\M'_{k_0}$
is concerned, we may replace $c_0$ by $c'_0$. 

We will now use internal addresses \cite{IntAddr}. The internal
address of any hyperbolic component is finite by definition (it is a
strictly increasing sequence of integers starting with $1$ and ending
with the period of the component). Any hyperbolic component
containing $c'_0$ in its wake but not within its subwake at internal
angle $1/2$ must appear in the internal address of $c'_0$ by
\cite[Lemma~6.4]{IntAddr}. This can happen only finitely often. (The
reason for this is that upon entering the $p/q$-subwake of any such
hyperbolic component, a periodic orbit is created, and every point of
this orbit is the landing point of $q$ dynamic rays. Every hyperbolic
Julia set has only finitely many such orbits with $q>2$). All but
finitely many little Mandelbrot sets $\M'_k$ will thus contain $c'_0$
and $c_0$ within their sublimbs at internal angles $1/2$.

Let $c_{-2}$ be the main antenna tip of $\M'_{k_0}$. We apply the
Branch Theorem \cite[\ThmBranch]{FiberMulti} to $c_{-2}$ and $c'_0$.
These two points are either separated by a hyperbolic component or by
a Misiurewicz point which must be on the tuned image of the real axis
within $\M'_{k_0}$. However, Misiurewicz points on the interior of
the tuned image of the real axis are never branch points by
Lemma~\ref{LemBranchReal}, so the separation takes place at a
hyperbolic component. This hyperbolic component must then show up in
the internal address of $c'_0$ for the same reason as above. It
follows that the third case above can happen only finitely many
times. Ignoring finitely many $\M'_k$, we can completely ignore the
third case.

If at least one of the $\M'_k$ realizes the second possibility, we
are done: the tuning map of $\M'_k$ is a homeomorphism from $\M$ onto
$\M'_k$ and we can use its restriction to the real axis of the
Mandelbrot set. This connects the root to the main antenna tip of
$\M'_k$, and composition with the projection will connect the root
to the main antenna tip of the projection of $\M'_k$. Therefore,
$\phi$ is continuous in a neighborhood of $t^*$. By assumption on
the map $\phi$, the composition $\pi\circ\psi$ coincides with $\phi$
(possibly up to reparametrization) within the entire neighborhood of
$t^*$ in the tuned image of the real axis. 

Otherwise, the first case must happen infinitely often. It follows
that all $\M'_k$ must contain $c_0$ as well as $c^*$, so both points
are in the same fiber. Since we assumed that the fiber of $c_0$ is
trivial, we have $\{c_0\}=\pi^{-1}(c^*)$. But we can then simply
define $\psi(t^*):=c_0$.

This way, we have defined $\psi\colon[0,1]\to\M$, and this map is
continuous by construction. This finishes the proof.
\qed

This result seems to have been observed first by Kahn, at least for
the case of dyadic Misiurewicz points. A sketch of proof, from which
the idea of this proof is taken, can be found in 
Douady~\cite{DoCompacts}.

\hide{
It seems quite possible that an analogous theorem is true for all
Multibrot sets. While it might be difficult to extend Yoccoz' puzzle
arguments, the approach by Branner and Douady might go through.
}

For the real line of the Mandelbrot set, there are many results
known, most of which are quite recent: hyperbolicity is dense
(Lyubich~\cite{LyHypDense} and Graczyk and Swi\c{a}tek \cite{GS}),
and all the Julia sets are locally connected (Levin and van Strien
\cite{LvS}, with simplifications by Sands, and Lyubich and Yampolsky
\cite{LY}). Therefore, we can draw the following conclusions:

\begin{corollary}[Properties of Paths in Mandelbrot Set]
\label{CorPropertiesMandelPath} \lineclear
The paths in the Mandelbrot set, as defined in the previous theorem,
run only through parameters with locally connected Julia sets, and
hyperbolicity is dense on these paths.
\qedd
\end{corollary}

The Branch Theorem (\cite[Theorem~9.1]{IntAddr} or
\cite[\ThmBranch]{FiberMulti}) shows that the various paths within
the Mandelbrot set starting at the origin and leading to different
points can split only at Misiurewicz points or within hyperbolic
components. 

There is a concept of ``regular arcs'' or ``legal arcs'', due to
Douady and Hubbard~\cite{Orsay} (see also \cite{DoCompacts}): a
regular arc is an arc within the Mandelbrot set subject to the
condition that it traverses any hyperbolic component only along the
union of two internal rays: it may enter the component along an
internal ray towards the center, and then leave the component along
another internal ray. It is easy to see that a regular arc in the
Mandelbrot set is uniquely specified by its endpoints (except
possibly for some freedom within non-hyperbolic components --- but
that will never happen for the arcs constructed here because
hyperbolicity is dense on them). Branch points of regular arcs are
necessarily postcritically finite.

Here are two more corollaries about the Mandelbrot set.
\begin{corollary}[Crossed Renormalization and Trivial Fibers]
\label{CorCrossedTunFibers} \lineclear
The fiber of any crossed renormalizable parameter in the Mandelbrot
set is trivial if and only if the fiber of the crossed renormalized
parameter is trivial.
\end{corollary}
\proof
According to \cite[Section~3]{CrossRenorm}, the subset of the
Mandelbrot set corresponding to any particular type of crossed
renormalization is homeomorphic to a sublimb of the Mandelbrot set,
and it is separated from the rest of the Mandelbrot set by a
countable collection of rational parameter ray pairs, just like in
Theorem~\ref{ThmLittleMandelBoundary}. Therefore, the fiber of any
crossed renormalizable parameter is contained within the crossed
renormalization locus of the same type. Within this locus, the same
argument as in Corollary~\ref{CorFiberTuning} applies.
\qed

\begin{corollary}[Branner-Douady Homeomorphisms and Trivial Fibers]
\label{CorBDFibers} \lineclear
The Branner-Douady homeomorphism from the $1/2$-limb into the
$1/3$-limb of the Mandelbrot set preserves triviality of fibers of\/
$\M$.
\end{corollary}
\sketch
The proof is essentially the same as above, so we just sketch the
main steps. The first step is again to show that the image of the
$1/2$-limb within the $1/3$-limb is bounded by rational parameter ray
pairs, and this follows from the corresponding properties in the
dynamic planes by the same transfer arguments as for simple and
crossed renormalization.

We can then conclude that the fiber of any point in the image of the
homeomorphism is entirely contained within the image. Finally, we
show as before that any separation line in the $1/2$-limb carries
over to a separation line in the $1/3$-limb with the corresponding
separation properties in the image, and conversely.
\qed

\remark
Another way to say this is like this: the Branner-Douady
homeomorphism is compatible with tuning, so that it maps little
Mandelbrot sets to little Mandelbrot sets. Points which are not within
those have trivial fibers by Yoccoz' result, and for points within
little Mandelbrot sets, the Branner-Douady image is the same as the
image under tuning maps. Extra decorations are attached within the
$1/3$-limb to the image of the $1/2$-limb. These are separated from
the little Mandelbrot set by rational parameter ray pairs, so fibers
cannot get larger under this homeomorphism. However, this argument
requires the Yoccoz theorem and is therefore restricted to degree two
only.

\medskip

Finally, we give a result about Siegel disks. The usual choice of
rational external angles in the definition of fibers ($Q=\Q/\Z$) had
to be extended by the grand orbit of the ray landing at the critical
value precisely when there is a Siegel disk. We will now give a use
of this extra work. First we recall a recent theorem of
Petersen~\cite{Pt}, reproved by Yampolsky. A Siegel disk is called of
{\em bounded type} if its multiplier $\mu=e^{2\pi i \theta}$ is such
that the continued fraction expansion of $\theta$ has bounded entries.

\begin{theorem}[Local Connectivity of Period One Siegel Disks]
\label{ThmCarstenSiegel} \lineclear
For every quadratic polynomial with a period one Siegel disk of
bounded type, the Julia set is locally connected.
\end{theorem}

We can now remove the condition on the period.
\begin{corollary}[Trivial Fibers and Bounded Type Siegel Disks]
\label{CorSiegelGeneral} \lineclear
For every quadratic polynomial with a Siegel disk of bounded type,
the Julia set has trivial fibers and is thus locally connected.
\end{corollary}
\proof
The Siegel disk of period one is locally connected and thus has
trivial fibers for an appropriate choice of $Q$
\cite[\PropLocConnJuliaFiber]{Fibers}. Any quadratic Siegel disk of
any period $n\geq 1$ has its parameter on the boundary of a
hyperbolic component of period $n$, so it is the image of a parameter
on the boundary of the main cardioid of the Mandelbrot set under a
tuning map of period $n$. The dynamics is thus $n$-renormalizable.
The result now follows from Corollary~\ref{CorFiberRenormJulia}.
\qed

\bigskip
\noindent Dierk Schleicher\\
\medskip
Fakult\"at f\"ur Mathematik\\ Technische
Universit\"at\\ Barer Stra{\ss}e 23\\ D-80290 M\"unchen, Germany\\
{\sl dierk$@$mathematik.tu-muenchen.de} 


\end{document}

%% file: LaTeX-top.tex
\newfont{\bbf}{msbm10 at 12pt}

\def\Q{\mbox{\bbf Q}}
\def\Z{\mbox{\bbf Z}}
\def\R{\mbox{\bbf R}}
\def\C{\mbox{\bbf C}}
\def\Cbar{\overline{\cz}}

\def\disk{\mbox{\,\bbf D}}

\def\diskbar{\ovl{\disk}}

\def\cz{\mbox{\bbf C}}

\def\ovl{\overline}
\def\phi1{\phi}
\def\phi{\varphi}

\def\*{ {\mbox{\tt $\star$}} }
\def\0{ {\mbox{\tt 0}} }
\def\1{ {\mbox{\tt 1}} }
\def\2{ {\mbox{\tt 2}} }
\def\3{ {\mbox{\tt 3}} }
\def\4{ {\mbox{\tt 4}} }

\def\St#1#2 {{\small\tt\rule{0pt}{0pt}_{\mbox{#1}}^{\mbox{#2}} }}

\def\IntAdr{\rightarrow}

\def\<{\prec}
\def\>{\succ}

\ifx\ThmNum\undefined 
   \newtheorem{theorem}{Theorem}[section]
   \def\newsection#1{\setcounter{theorem}{0} \section{#1}}
\else 
   \newtheorem{theorem}{Theorem} 
   \def\newsection#1{\section{#1}}
\fi

\newtheorem{proposition}[theorem]{Proposition}
\newtheorem{lemma}[theorem]{Lemma}
\newtheorem{definition}[theorem]{Definition}

\newtheorem{corollary}[theorem]{Corollary}

\def\proof{\par\medskip\noindent {\sc Proof. }}
\def\proofof #1 {\par\medskip\noindent {\sc Proof of #1. }}
\def\sketch{\par\medskip\noindent{\sc Sketch of proof. }}
\def\sketchof #1 {\par\medskip\noindent {\sc Sketch of proof of #1. }}
\def\Box{\framebox[10pt]{\rule{0pt}{3pt}}}
\def\qed{\hfill $\Box$ \medskip \par}
\def\qedd{\rule{0pt}{0pt}\hfill $\Box$}
\def\remark{\par\medskip \noindent {\sc Remark. }}
\def\lineclear{\rule{0pt}{0pt}\par\noindent}

\def\LabelCaption#1#2{
   \centerline{\parbox{\captionwidth}{   
   \caption{\sl #2}  \label{#1} } }
}

\def\reminder #1 {{\sf #1}}
\def\hide #1 {}
\long\def\longhide #1 {}

%% file: imsmark.tex
\def\IMSmarkvadjust{0 pt}
\def\IMSmarkhadjust{0 pt}
\def\SBIMSMark#1#2#3{
 \font\SBF=cmss10 at 10 true pt
 \font\SBI=cmssi10 at 10 true pt
 \setbox0=\hbox{\SBF Stony Brook IMS Preprint \##1}
 \setbox2=\hbox to \wd0{\hfil \SBI #2}
 \setbox4=\hbox to \wd0{\hfil \SBI #3}
 \setbox6=\hbox to \wd0{\hss
             \vbox{\hsize=\wd0 \parskip=0pt \baselineskip=10 true pt
                   \copy0 \break%
                   \copy2 \break%
                   \copy4 \break}}
 \dimen0=\ht6   \advance\dimen0 by \vsize \advance\dimen0 by 8 true pt
                \advance\dimen0 by -\pagetotal
	        \advance\dimen0 by \IMSmarkvadjust
 \dimen2=\hsize \advance\dimen2 by .25 true in
	        \advance\dimen2 by \IMSmarkhadjust

%
%
  \openin2=publishd.tex
  \ifeof2\setbox0=\hbox to 0pt{}
  \else 
     \setbox0=\hbox to 3.1 true in{
                \vbox to \ht6{\hsize=3 true in \parskip=0pt  \noindent  
                {\SBI Published in modified form:}\hfil\break
                \input publishd.tex 
                \vfill}}
  \fi
  \closein2
  \ht0=0pt \dp0=0pt
 \ht6=0pt \dp6=0pt
 \setbox8=\vbox to \dimen0{\vfill \hbox to \dimen2{\copy0 \hss \copy6}}
 \ht8=0pt \dp8=0pt \wd8=0pt
 \copy8
 \message{*** Stony Brook IMS Preprint #1, #2. #3 ***}
}